\documentclass{article}

\usepackage{amsthm,amsmath,amsfonts,amssymb}
\usepackage{graphicx}
\usepackage{hyperref}

\newtheorem{theorem}{Theorem}[section]
\newtheorem{lemma}[theorem]{Lemma}
\newtheorem{corollary}[theorem]{Corollary}
\newtheorem{remark}[theorem]{Remark}
\newtheorem{definition}[theorem]{Definition}
\newtheorem{proposition}[theorem]{Proposition}
\newtheorem{example}[theorem]{Example}

\newcommand{\R}{{\mathbb R}}

\def \ts {\textstyle}
\def \eps {\epsilon}
\def \ds {\displaystyle}
\def \half{\frac{1}{2}}
\newcommand{\dwto}{{\stackrel{E}{-{\hspace{-2mm}}\rightharpoonup}}}
\newcommand{\dto}{{\stackrel{E}{\longrightarrow}}}

\begin{document}

\title{Concentrated terms and varying domains in elliptic equations: Lipschitz case}

\author{Gleiciane S.Arag\~ao\thanks{Departamento de Ci\^encias Exatas e da Terra, Universidade Federal de S\~ao Paulo. Rua Professor Artur Riedel, 275, Jardim Eldorado, Cep 09972-270, Diadema-SP, Brazil, Phone: (+55\ 11) 33193300. E-mail: gleiciane.aragao@unifesp.br. Partially supported by CNPq 475146/2013-1, Brazil.} \, and\,  Simone M. Bruschi\thanks{Departamento de Matem{\'a}tica, Universidade de Bras\'ilia. Campus Universit\'ario Darcy Ribeiro, ICC Centro, Bloco A, Asa Norte, Cep 70910-900, Bras\'ilia-DF, Brazil, Phone: (+55\ 61) 31076389. E-mail: sbruschi@unb.br. Partially supported by FEMAT, Brazil.}}

\date{}
\maketitle \thispagestyle{empty} \vspace{-10pt}

\date{}
\maketitle \thispagestyle{empty} \vspace{-10pt}

\begin{abstract}
In this paper, we analyze the behavior of a family of solutions of a nonlinear elliptic equation with nonlinear boundary conditions, when the boundary of the domain presents a highly oscillatory behavior which is uniformly Lipschitz and  nonlinear terms are concentrated in a region which neighbors  the boundary  domain. We prove that this family of solutions converges to the solutions of a limit problem in $H^{1}$, an elliptic equation with nonlinear boundary conditions which  captures the oscillatory behavior of the boundary and whose nonlinear terms are transformed into a flux condition on the boundary. Indeed, we show the upper semicontinuity of this family of solutions.

\vspace{0.2cm}


\vspace{0.2cm}

\end{abstract}

\maketitle

\vspace{-6pt}


\section{Introduction}
\label{introd}

\par In this work, we analyze the behavior of the family of solutions of a concentrated elliptic equation with nonlinear boundary conditions of the type
\begin{equation}  
\label{nbc}
\left\{
\begin{array}{ll}
-\Delta u_{\epsilon}+u_{\epsilon}=\displaystyle \frac{1}{\epsilon} \mathcal{X}_{\omega_{\epsilon}}f(x,u_{\epsilon}), &\hbox{in} \ \Omega_\epsilon \\
\frac{\ts\partial u_{\epsilon}}{\ts\partial n}+g(x,u_{\epsilon})=0, & \hbox{on} \
\partial \Omega_\epsilon
\end{array} \right.
\end{equation}
when the boundary of the domain presents a highly oscillatory behavior, as the parameter $\epsilon\to 0$, and nonlinear terms are concentrated in a region of the domain neighboring the boundary. To describe the problem, we will consider a family of uniformly bounded smooth domains $\Omega_\epsilon \subset \R^N$, with $N\geq 2$ and $0\leq\eps\leq \epsilon_0$, for some $\epsilon_0>0$ fixed, and we will look at this problem from the perturbation of the domain point of view and we will refer to $\Omega \equiv\Omega_0$ as the unperturbed domain and $\Omega_\epsilon$ as the perturbed domains. We will assume that $\Omega \subset \Omega_{\epsilon}$, $\Omega_\epsilon\to \Omega$ and $\partial\Omega_\epsilon\to \partial\Omega$ in the sense of Hausdorff, that is, $\mbox{dist}(\Omega_\epsilon,\Omega)+
\mbox{dist}(\partial\Omega_\epsilon,\partial\Omega)\to 0$ as $\epsilon\to 0$, where $\mbox{dist}$ is the symmetric Hausdorff distance of two sets in $\R^N$ ($\mbox{dist}(A,B)= \sup_{x\in A}\inf_{y\in B}|x-y|+\sup_{y\in B}\inf_{x\in A}|x-y|$). We will also assume that the nonlinearities $f,g:U\times\R\to\R$ are continuous in both variables and $C^2$ in the second one, where $U$ is a fixed and smooth bounded domain containing all $\overline{\Omega}_\epsilon$, for all $0\leq\eps\leq\eps_0$. Now, for sufficiently small $\epsilon$, $\omega_\epsilon$ is the region between the boundaries of $\partial \Omega$ and $\partial \Omega_{\epsilon}$. Note that $\omega_{\epsilon}$ shrinks to $\partial \Omega$ as $\epsilon\to 0$ and we use the characteristic function $\mathcal{X}_{\omega_{\epsilon}}$ of the region $\omega_{\epsilon}$ to express the concentration in $\omega_{\epsilon}$. The Figure \ref{figure1} illustrates the oscillating set $\omega_{\epsilon} \subset \overline{\Omega}_{\epsilon}$.

\begin{figure}[h]
\label{figure1}
\begin{center}
\includegraphics[width=0.38\linewidth,height=0.36\textheight,keepaspectratio]{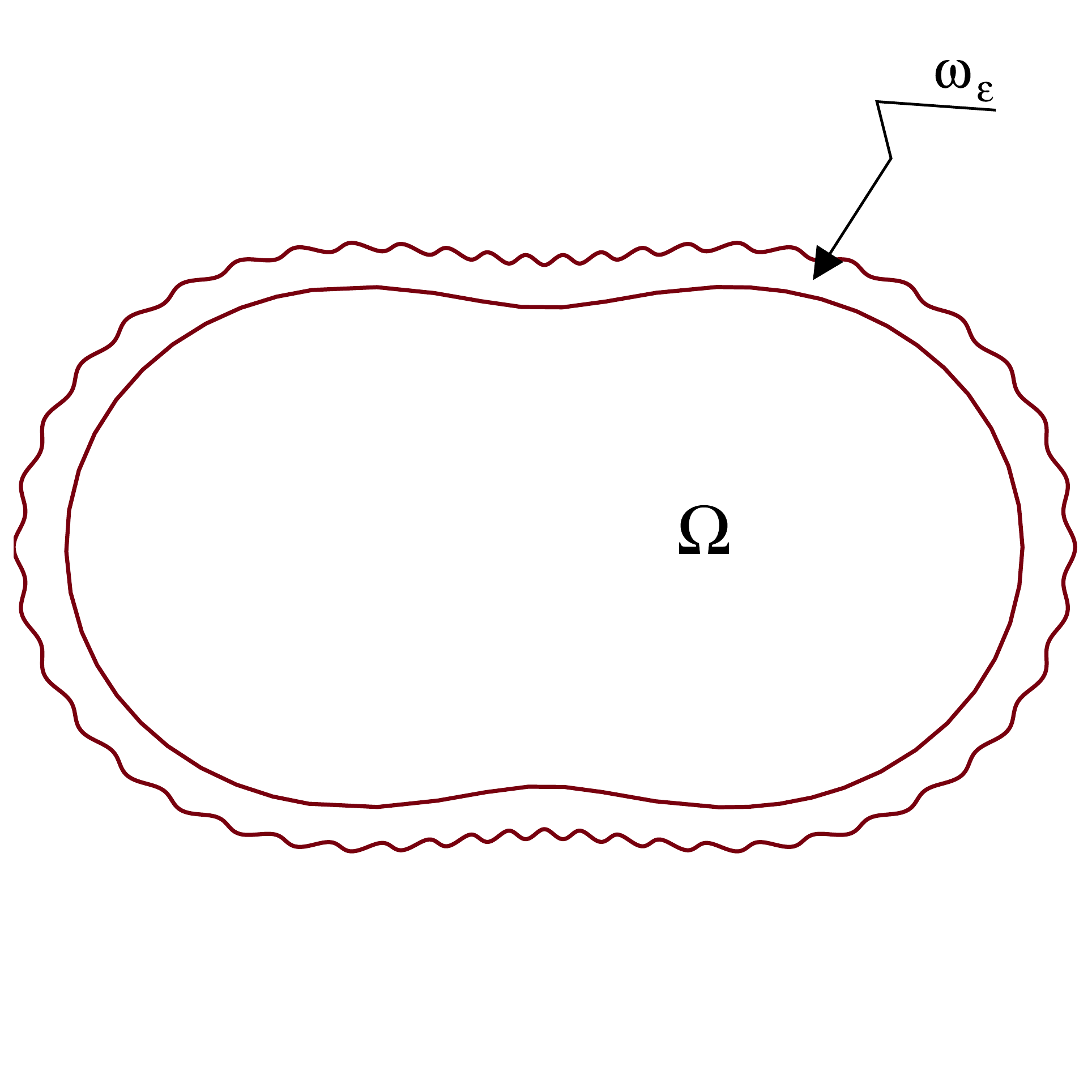}
\vspace*{-1.2cm}
\caption{The set $\omega_{\epsilon}$.}
\label{figure1}
\end{center}
\end{figure}

Although the domains behave continuously  as $\epsilon\to 0$, the same is not true for the boundaries. It is possible that the measure $|\partial\Omega_\epsilon|$ is uniformly bounded and that it  does not converge to $|\partial\Omega|$ when $\epsilon \to 0$. In fact, this is the case when  the boundary $\partial\Omega_\epsilon$ presents an oscillatory behavior in which, up to a diffeomorphism, the period goes to zero with the same order of  amplitude. 

Since $\omega_{\epsilon}$ shrinks to $\partial \Omega$ as $\epsilon\to 0$,  it is reasonable to expect that the family of solutions $\{u_{\epsilon}\}$ of (\ref{nbc}) will converge to a solution of an equation  with a nonlinear boundary condition on $\partial \Omega$ that  inherits the information about the  concentration' region. Moreover, the boundary condition in the limit problem also inherits the information about the behavior  of the measure of the deformation of $\partial \Omega_\epsilon$ with respect to $\partial \Omega$. We will show that the solutions of (\ref{nbc}) converge in $H^1(\Omega_{\epsilon})$ to the solutions of the elliptic problem 
\begin{equation}  
\label{nbc_limite_gamma_F}
\left\{
\begin{array}{ll}
-\Delta u+u=0, &\hbox{in} \ \Omega \\
\frac{\ts\partial u}{\ts\partial n}+\gamma(x)g(x,u)=\beta(x)f(x,u), & \hbox{on} \
\partial \Omega
\end{array} \right.
\end{equation}
where the function $\gamma\in L^\infty(\partial\Omega)$  is related to the behavior of the measure $(N-1)$-dimensional of the $\partial\Omega_\epsilon$  and $\beta \in L^\infty(\partial\Omega)$ is related to the behavior of the  measure $N$-dimensional  of   the  concentration' region $\omega_\epsilon$. Indeed, we will prove the upper semicontinuity of the family of solutions of (\ref{nbc}) and (\ref{nbc_limite_gamma_F}) in $H^1(\Omega_\epsilon)$. The precise hypotheses about  the domains, the exact definition of  $\omega_{\epsilon}$ and the functions $\gamma$ and $\beta$ are given in Section \ref{perdominio}.  

The behavior of the solutions of nonlinear elliptic equations with nonlinear boundary conditions and rapidly  varying  boundaries  was studied in \cite{AB0}, for the case of uniformly  Lipschitz deformation of the boundary, but  without concentration. The authors proved the convergence of the solutions in $H^1(\Omega_{\epsilon})$ and $C^{\alpha}(\Omega_{\epsilon})$, for some $0<\alpha\leq 1$.

The behavior of the solutions of elliptic problems with terms concentrated in a neighborhood of the boundary of the  domain  was initially studied in~\cite{arrieta}, when the neighborhood is a strip of width $\epsilon$ and has a base in the boundary, without oscillatory behavior and inside of $\Omega$. More precisely, in~\cite{arrieta} it was proven that these solutions converge, in certain Bessel Potentials spaces and in the continuous functions space, to the solution of an elliptic problem where the reaction term and the concentrating potential are transformed into a flux condition and a potential on the boundary. Later, the asymptotic behavior of the attractors of a parabolic problem was analyzed in~\cite{anibal}, where the upper semicontinuity of attractors was proved. In  \cite{arrieta,anibal} the domain $\Omega$ is $C^2$ in $\R^N$.

Recently, in~\cite{gam1} some results of~\cite{arrieta} were adapted to a nonlinear elliptic problem posed on an open square $\Omega$ in $\R^2$, considering $\omega_\epsilon \subset \Omega $ and with  highly oscillatory behavior in the boundary inside of $\Omega$. The dynamics of the flow generated  by a nonlinear parabolic problem posed on a $C^2$ domain $\Omega$ in $\R^2$, when some reaction and potential terms are concentrated in a neighborhood of the boundary and the ``inner boundary'' of this neighborhood presents a highly oscillatory behavior, was studied in \cite{gam2} where the continuity of the family of attractors was proved. 

It is important to note that all previous works with terms concentrating in a neighborhood of the boundary deal with non varying domain and since $\omega_\epsilon$ is inside of $\Omega$ then all the equations are defined in the same domain. In our case, the region $\omega_\epsilon$  is outside of $\Omega$. We also generalize the domain treated in the literature allowing  $\Omega, \Omega_\epsilon \subset \R^N$, with $N\geq 2$,  and alowing $\Omega$ to be a $C^1$ domain. Moreover, the oscillations discussed in this work are more general than those in~\cite{gam1,gam2}.

This paper is organized as follow: in Section \ref{perdominio}, we define the domain perturbation and we state our main result (Theorem \ref{th-main}). In Section \ref{technical-results}, we obtain the solutions of (\ref{nbc}) and (\ref{nbc_limite_gamma_F}) as fixed points of appropriate nonlinear maps. We also obtain some technical results on uniform estimates in Bessel Potentials spaces and of constant of trace operators, and we analyze the limit of concentrated integrals. In Section \ref{upper}, we  prove the upper semicontinuity of the family of solutions of (\ref{nbc}) and (\ref{nbc_limite_gamma_F}). In Section \ref{conclusion}, we state some results that can be proven, but we  leave out  the proof here.


\section{Setting of the problem and main results}
\label{perdominio}

We consider a family of smooth bounded domains $\Omega_\epsilon  \subset \R^N$, with $N\geq 2$ and $0 \leq \epsilon  \leq \epsilon_0$, for some $\epsilon_0 >0$ fixed, and we regard $\Omega_\epsilon$ as a perturbation of the fixed domain $\Omega\equiv\Omega_0$. We consider the following hypothesis on the domains
\begin{description}
\item[(P)] There exists a finite open cover $\{U_i\}_{i=0}^m$ of $\Omega$ such that $\overline{U}_0 \subset \Omega$, $\partial\Omega \subset \cup_{i=1}^m U_i$ and for each $i=1, \ldots, m$, there exists a Lipschitz diffeomorphism $\Phi_{i}: Q_N \to U_i$, where $Q_N=(-1,1)^N\subset \R^N$, such that
$$
\Phi_{i}(Q_{N-1} \times (-1,0)) = U_i \cap \Omega \qquad \mbox{and} \qquad \Phi_{i}(Q_{N-1} \times \{0\}) = U_i \cap \partial \Omega.
$$
We assume that $\overline{\Omega}_\epsilon \subset \cup_{i=0}^m U_i\equiv U$. For each $i=1, \ldots, m$, there exists a Lipschitz  function  $\rho_{i,\epsilon}: Q_{N-1} \to (-1,1)$ such that $\rho_{i,\epsilon} \to 0$ as $\epsilon \to 0$, uniformly in $Q_{N-1}$, and $\| \nabla \rho_{i,\epsilon}\|_{L^\infty(Q_{N-1})}\leq C$, with $C>0$ independent of $\epsilon$, $i=1,\ldots,m$.

Moreover, we assume that $\Phi_{i}^{-1}(U_i \cap \partial \Omega_\epsilon)$ is the graph of $\rho_{i,\epsilon}$ this means 
$$
U_i \cap \partial \Omega_\epsilon=\Phi_{i}(\{(x',\rho_{i,\epsilon}(x')) \;:\; x' \in  Q_{N-1}\}).
$$
\end{description}

\begin{figure}[h]
\begin{center}
\includegraphics[width=0.38\linewidth,height=0.36\textheight,keepaspectratio]{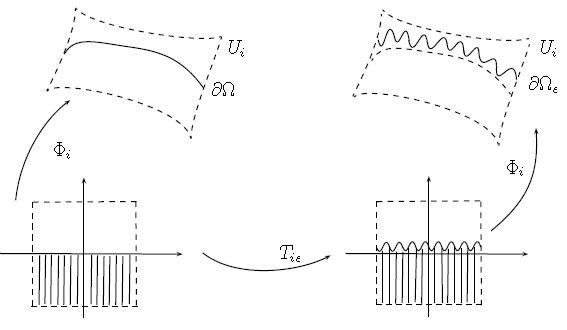}
\vspace*{0cm}
\caption{The parametrizations.}
\label{figure2}
\end{center}
\end{figure}

We consider the following  mappings $T_{i,\epsilon}:Q_N \to Q_N$ defined by
$$
T_{i,\epsilon}(x',s) =
\left\{
\begin{array}{ll}
(x',s+s\rho_{i,\epsilon}(x')+\rho_{i,\epsilon}(x')),&\hbox{ for }s \in (-1,0) \\
(x',s-s\rho_{i,\epsilon}(x')+\rho_{i,\epsilon}(x')),&\hbox{ for }s \in [0,1).
\end{array}
\right.
$$
Also, 
$$
\Phi_{i,\epsilon}:=\Phi_{i}\circ T_{i,\epsilon}:Q_N \to U_i 
$$
and we also denote by 
$$
\begin{array}{rl}
\psi_{i,\epsilon}:Q_{N-1}&\to U_i\cap\partial\Omega_\epsilon \\
x'& \mapsto \Phi_{i,\epsilon}(x',0)
\end{array}
\qquad \mbox{and} \qquad \begin{array}{rl}
\psi_{i}:Q_{N-1}&\to U_i\cap\partial\Omega \\
x'& \mapsto \Phi_{i}(x',0).
\end{array}
$$
Notice that $\psi_{i,\epsilon}$ and $\psi_{i}$ are local parametrizations of $\partial\Omega_\epsilon$ and $\partial\Omega$, respectively. Furthermore, observe that all the maps above are Lipschitz. The Figure \ref{figure2} illustrates the parametrizations.

With the notations above, we define 
$$
\omega_{\epsilon}=\bigcup_{i=1}^{m}\Phi_{i}\left( \left\{(x',x_N)\in \mathbb{R}^N \;:\; 0\leq x_N < \rho_{i,\epsilon}(x') \quad \mbox{and}\quad x'=(x_1,...,x_{N-1})\in Q_{N-1} \right\} \right), \quad \mbox{for $0<\epsilon \leq \epsilon_0$}.
$$

In some results, specially in the main result,  we consider additional  hypotheses on the deformation $\partial\Omega_\epsilon$. 
In order to give the hypothesis to deal with the concentration we need the following definition

\begin{definition}
\label{definition-jacobian}
Let $\eta: A \subset \R^{N-1} \to \R^N$ almost everywhere differentiable,  we define  
the $(N-1)$-dimensional Jacobian of $\eta$  as
$$
J_{N-1}\eta \equiv \left|{\partial \eta \over \partial x_1}\wedge\ldots\wedge {\partial \eta\over \partial x_{N-1}}\right| = \sqrt{\sum_{j=1}^N (\mbox{det}({\rm{Jac}}\,\eta)_j)^2},
$$
where $v_1 \wedge\ldots\wedge v_{N-1}$ is the exterior product of the $(N-1)$ vectors $v_1,\ldots,v_{N-1} \in \mathbb{R}^{N}$
and $({\rm{Jac}}\,\eta)_j$ is the $(N-1)$-dimensional matrix obtained by deleting the $j$-th row of the Jacobian matrix of $\eta$.  
\end{definition}

We use $J_N$ for the absolute value of the $N$-dimensional Jacobian determinant. Now, we are in condition to give the hypothesis.

\begin{description}  
\item[(B)] For each $i=1, \ldots, m$, $\rho_{i,\epsilon}$ is $O(\epsilon)$ when $\epsilon \to 0$
and  there exists  a function $\tilde{\beta}_i \in L^\infty(Q_{N-1})$ such that 
$$
\frac{1}{\epsilon}\rho_{i,\epsilon} \stackrel{\epsilon \to 0}{-{\hspace{-2mm}}\rightharpoonup}    \tilde{\beta}_i \quad \mbox{in $L^1(Q_{N-1})$}.
$$
\end{description}

\begin{definition}
\label{definition-of-beta}
For  $x\in U_i\cap\partial\Omega$, let $(x^\prime,0) =\Phi_i^{-1} (x) \in Q_{N}$, we define $\beta:\partial\Omega \to \R$ as 
$$
\beta(x)=\frac{\tilde{\beta}_i(x^\prime)(J_N\Phi_i)(x^\prime,0)}{J_{N-1}\psi_i(x^\prime)}.
$$
\end{definition}

We observe that, since $J_N(\Phi_i(x^\prime, \rho_{i, \epsilon}(x^\prime)z)) = \rho_{i, \epsilon} (x^\prime)(J_N\Phi_i) (x^\prime, \rho_{i, \epsilon}(x^\prime)z)$, $z\in [0,1]$, the hypothesis {\bf(B)} means that 
$$
\frac{1}{\epsilon} \frac{J_N(\Phi_i(x^\prime, \rho_{i, \epsilon}(x^\prime)z))}{J_{N-1}\psi_i(x^\prime)}  \stackrel{\epsilon \to 0}{-{\hspace{-2mm}}\rightharpoonup} \beta(x).
$$

We note that $\beta$ could depend on the choice of the charts $U_i$ and the maps $\Phi_i$ and $\rho_{i,\epsilon}$. We will prove in the Corollary \ref{inde2} that $\beta$ is well defined and unique for the family $\Omega_\epsilon$, $0\leq \epsilon \leq \epsilon_0$.

We give an example of the function $\rho_{i,\epsilon}$ satisfying the hypothesis {\bf(B)}  and its correspondent function $\tilde{\beta}_i$. 
\begin{example}
For each $i=1, \ldots, m$, let $\rho_{i}: \R^{N-1} \to \mathbb{R}^+$ be a $Y$-periodic Lipschitz function, where $Y=(0,l_1)\times \cdots \times (0, l_{N-1}) \in \R^{N-1}$ with $l_1,\ldots,l_{N-1} \in \R^{+}$ (a function $\rho_{i}$ is called $Y$-periodic iff $\rho_{i}(x'+kl_je_j)=\rho_{i}(x')$ on $Q_{N-1}$, $\forall k\in \mathbb{Z}$ and $\forall j\in \{1,\ldots,N-1\}$, where $x'=(x_1, \ldots, x_{N-1})\in \R^{N-1}$ and $\{e_1,\ldots,e_{N-1}\}$ is the canonical basis of $\mathbb{R}^{N-1}$), and we define $\rho_{i,\epsilon}: Q_{N-1} \to [0,1)$ by
$$
\rho_{i,\epsilon}(x')=\displaystyle \epsilon \varphi(x^\prime) \rho_{i}\left(\frac{x'}{\epsilon^\alpha}\right), 
$$
for $x'\in Q_{N-1}$, sufficiently small $\epsilon$, say $0<\epsilon \leq \epsilon_0$, $0<\alpha \leq 1$ and $\varphi: Q_{N-1} \to \R$ a continuous function. From Theorem 2.6 in \cite{cioranescu}, we get
$$
\frac{1}{\epsilon}\rho_{i,\epsilon} \stackrel{\epsilon \to 0}{-{\hspace{-2mm}}\rightharpoonup}  \varphi M_Y(\rho_i)=\tilde{\beta}_i \quad \mbox{in $L^1(Q_{N-1})$}, 
$$
where $M_{Y}(\rho_i)$ is the mean value of $\rho_i$ over $Y$ given by
$$
M_{Y}(\rho_i)=\frac{1}{\left|Y\right|}\int_{Y} \rho_i(x')dx'.
$$
 \end{example}

In order to deal with the interaction between the nonlinear boundary condition and the oscillatory behavior of $\partial\Omega_\epsilon$  we consider the hypothesis

\begin{description}
\item[(G)] For each $i=1,\ldots,m$, there exists a function $\tilde{\gamma}_i \in L^\infty(Q_{N-1})$ such that
$$J_{N-1}\psi_{i,\epsilon} \stackrel{\epsilon \to 0}{-{\hspace{-2mm}}\rightharpoonup} \tilde{\gamma}_i \quad \mbox{in $L^1(Q_{N-1})$}.$$
\end{description}

\vspace{0.2cm}

Considering the hypotheses {\bf(P)} and {\bf(G)} let  $\gamma:  \partial\Omega \to \R$ be a function which measure the limit of the deformation of $\partial \Omega_\epsilon$ relatively to $\partial\Omega$ (see \cite{AB0}). More precisely, we have

\begin{definition}
\label{definition-of-gamma}
For $x\in U_i\cap\partial\Omega$, let $(x',0) =\Phi_i^{-1} (x) \in Q_{N}$, we define $\gamma:  \partial\Omega \to \R$ as
$$
\gamma(x) = \frac{\tilde{\gamma}_i(x')}{J_{N-1}\psi_i(x')}.$$
\end{definition}

Also, the function $\gamma$ is independent of $\psi_{i,\epsilon}$, $U_i$ and the maps $\Phi_i$. This was proved in the Corollary 5.1 in \cite{AB0}. 

We  consider an example where we can  explicitly calculate the functions $\gamma$ and $\beta$.
\begin{example}
\label{basic-example}
Consider the domain $\Omega \subset \R^2$. Let  $\Phi: (-1,1)^2 \to \R^2$ given by 
$$
\Phi(x,y)= \left( \left(r(x)+\frac{y}{2} \right) \cos\left(\frac{\pi x}{2}\right),  \left(r(x)+\frac{y}{2} \right)\sin\left(\frac{\pi x}{2} \right) \right)
$$ be a chart of $\Omega$.  
Therefore, $J_2 \Phi(x,y) = \frac{\pi}{4}\left( r(x) +\frac{y}{2} \right)$ and $(J_2\Phi)(x,0)= \frac{\pi}{4}r(x)$. 

Since $\psi(x)=\Phi(x,0)= \left(r(x) \cos\left(\frac{\pi x}{2}\right),  r(x)\sin\left(\frac{\pi x}{2}\right) \right)$  then
$$J_1 \psi(x) =  \sqrt{(r^\prime(x))^2+ \frac{\pi^2}{4}(r(x))^2}.$$
 
Consider $\rho:\R \to \R^+$  $l$-periodic, let  $\rho_\epsilon(x)=\epsilon\rho(\frac{x}{\epsilon})$, 
$\Phi_\epsilon(x,y) = \Phi \circ T_\epsilon (x,y)$ and $\psi_\epsilon(x)= \Phi_\epsilon(x,0)$. We get, $
{\ds{ J_1\psi_\epsilon(x)=  \sqrt{\left(r^\prime(x)+\frac{\rho_\epsilon^\prime(x)}{2}\right)^2+\left(r(x)+\frac{\rho_\epsilon(x)}{2}\right)^2}}} $ and 
$$
{\ds{ J_1\psi_\epsilon(x)}}\stackrel{\epsilon \to 0}{-{\hspace{-2mm}}\rightharpoonup} \tilde{\gamma}(x)=  \frac{1}{l} \int_{0}^{l} \sqrt{\left(r^\prime(x) + \frac{\rho^\prime(z)}{2}\right)^2 + (r(x))^2} \; dz \quad \mbox{ in $L^1(-1,1)$}.
$$  
Then $\gamma(x_1,x_2) = \frac{\tilde{\gamma}(x)}{J_1\psi(x) }$ and $\beta(x_1,x_2)=\frac{M_l(\rho) (J_2\Phi)(x,0)}{J_1\psi(x)}$, where $(x,0)= \Phi^{-1}(x_1,x_2)$, $(x_1,x_2) \in \partial \Omega$. 

We observe that if $r(x)$ is constant then $\gamma$ and $\beta$ are  constants. This fact is expected since there is no distinction between the points on $\partial\Omega$ that is a circle.

Considering $\rho(x)=\sin (\pi x) + 2$ then 
$$
\tilde{\gamma}(x)=  \frac{1}{2\pi} \int_{0}^{2\pi} \sqrt{\left(r^\prime(x) + \frac{\pi\cos(z)}{2} \right)^2 + (r(x))^2} \; dz\quad  \mbox{in $L^1(-1,1)$}.
$$
\end{example}

Now, with respect to the equations, we will be interested in studying the behavior of the solutions of the elliptic equation (\ref{nbc}) where, as we mentioned in the introduction,  the nonlinearities $f, g:U\times\R\to\R$ are continuous in both variables and $C^2$ in the second one, where $U$ is a bounded domain containing $\overline{\Omega}_\epsilon$,  for all $0\leq\eps\leq\eps_0$.

Consider the family of spaces $H^1(\Omega_\epsilon)$ and $H^1(\Omega)$ with their usual norms. Since we will need to compare functions defined in $\Omega_\epsilon$ with functions defined in the unperturbed domain $\Omega\equiv\Omega_0$,  we will need a tool to compare functions which are defined in different spaces. The appropriate notion for this is the concept of $E$-convergence and a key ingredient for this will be the use of the extension operator $E_\epsilon:H^1(\Omega)\to H^1(\Omega_\epsilon)$, which is defined as $E_\epsilon = R_\epsilon \circ P$, where $P:H^1(\Omega) \to H^1(\R^N)$ is a linear and continuous operator that extends a function $u$ defined in $\Omega$ to a function defined in $\R^N$ and $R_\epsilon$ is the restriction operator from functions defined in $\R^N$ to functions defined in $\Omega_\epsilon$, $R_\epsilon w=w_{|\Omega_\epsilon}$. Observe that we also have $E_\epsilon:L^p(\Omega)\to L^p(\Omega_\epsilon)$ and  $E_\epsilon:W^{1,p}(\Omega)\to W^{1,p}(\Omega_\epsilon)$, for all $1\leq p\leq \infty$. Considering  $X_\epsilon=H^1(\Omega_\epsilon)$ or $L^p(\Omega_\epsilon)$ or $W^{1,p}(\Omega_\epsilon)$, for $\epsilon \geq 0$, from  \cite{AB0} we have
$$ 
\|E_\epsilon u\|_{X_\epsilon} \to \|u\|_{X_0} \mbox{ \ \ \ and \ \ \ } \|E_\epsilon\| \leq \|R_\epsilon\|\cdot \|P\|\leq \|P\|, \mbox{ \ \ \ independent of } \epsilon.
$$
The concept of $E$-convergence is defined as follows: $u_\epsilon {\stackrel{E}{\longrightarrow}} u$ if $\|u_\epsilon-E_\epsilon u\|_{H^1(\Omega_\epsilon)}\to 0$, as $\epsilon \to 0$.  We also have a notion of weak $E$-convergence, which is defined as follows:  $u_\epsilon {\stackrel{E}{-{\hspace{-2mm}}\rightharpoonup}} u$ if  $(u_\epsilon,w_\epsilon)_{H^1(\Omega_\epsilon)}\to (u,w)_{H^1(\Omega)}$, as $\epsilon \to 0$, for any sequence $w_\epsilon {\stackrel{E}{\longrightarrow}} w$, where $(\cdot,\cdot)_{H^1(\Omega_{\epsilon})}$ and  $(\cdot,\cdot)_{H^1(\Omega)}$ denote the inner product in $H^1(\Omega_{\epsilon})$ and $H^1(\Omega)$, respectively. More details about $E$-convergence can be found in Subsection 3.2 in \cite{AB0}.

Our main result is stated in the following theorem

\begin{theorem}
\label{th-main}
Assume {\bf (P)}, {\bf(B)} and {\bf (G)} are satisfied. Let $\{u_\epsilon^*\}$, $0<\epsilon\leq \epsilon_0$, be a family of solutions of (\ref{nbc}) satisfying $\| u_\epsilon^*\|_{L^\infty(\Omega_\epsilon)}\leq R$, for some constant $R>0$ independent of $\epsilon$. Then, there exist a subsequence $\{u_{\epsilon_k}^*\}$ and a function $u_0^* \in H^{1}(\Omega)$, with $\|u_0^*\|_{L^\infty(\Omega)}\leq R$, solution of (\ref{nbc_limite_gamma_F}) satisfying  $u^{*}_{\epsilon_k} {\stackrel{E}{\longrightarrow}} u^{*}_{0}$.
\end{theorem}

\begin{remark}
\label{rm-cut-off}
Since in Theorem \ref{th-main} we are concerned with solutions that are uniformly bounded in $L^\infty(\Omega_\epsilon)$, we may perform a cut-off in the nonlinearities $f$ and $g$ outside the region $|u|\leq R$, without modifying any of these solutions, in such a way that
\begin{equation}
\label{hip-f}
\begin{array}{l}
|f(x,u)|+|\partial_u f(x,u)|+|\partial_{uu}f(x,u)|\leq C \quad \mbox{and}\\ |g(x,u)|+|\partial_u g(x,u)|+|\partial_{uu}g(x,u)|\leq C,\quad \mbox{$\forall x\in U$ \; and \; $\forall u\in \R$}.
\end{array}
\end{equation}
\end{remark}


\section{Abstract setting and technical results}
\label{technical-results}

The solutions of (\ref{nbc}) and (\ref{nbc_limite_gamma_F}) will be obtained as fixed points of appropriate nonlinear maps defined in the spaces $H^1(\Omega_\epsilon)$ and $H^1(\Omega)$, respectively. These maps are constructed in Subsection \ref{fixed-points}. In the Subsections \ref{sobolev-trace} and \ref{results} we will also obtain several important technical results that will be needed in the proof of the main result. 


\subsection{Solutions as fixed points}
\label{fixed-points}

For $0<\epsilon\leq\eps_0$, consider the linear operator $A_\epsilon: D(A_\epsilon) \subset L^2(\Omega_\epsilon) \to
L^2(\Omega_\epsilon)$ defined by $A_\epsilon u_\epsilon = -\Delta u_\epsilon + u_\epsilon$ with domain $D(A_\epsilon) = \{u_\epsilon \in H^2(\Omega_\epsilon):
\frac{\partial u_\epsilon}{\partial n}=0  \mbox{\ in\ } \partial\Omega_\epsilon\}$. Let us denote by $X^0_\epsilon = L^2(\Omega_\epsilon)$, $X^1_\epsilon = D(A_\epsilon)$ and consider the scale of Hilbert spaces $\{(X^\alpha_\epsilon, A_{\epsilon,\alpha}) \; : \; \alpha \in \R\}$ constructed by complex interpolation, see \cite{Amann1}, which coincide, since we are in a Hilbert setting,  with the standard fractional power spaces of the operator $A_\epsilon$ and $X^{\alpha}_\epsilon \hookrightarrow H^{2\alpha}(\Omega_\epsilon)$, where $H^{s}(\Omega_{\epsilon})$ denotes the Bessel Potentials spaces of order $s\in \mathbb{R}$. This scale can also be extended to spaces of negative exponents by taking $X^{-\alpha}_\epsilon = (X^{\alpha}_\epsilon)^\prime$, for $\alpha >0$, and $X^{-\half}_\epsilon = H^{-1}(\Omega_\epsilon)$. Considering the realizations of $A_\epsilon$ in this scale, the operator $A_{\epsilon,-\half} \in \mathcal{L}(X^{\half}_\epsilon,X^{-\half}_\epsilon)$ is given by 
$$
\langle  A_{\epsilon,-\half}u_\epsilon, \phi_\epsilon\rangle = \int_{\Omega_\epsilon} (\nabla u_\epsilon \nabla \phi_\epsilon + u_\epsilon \phi_\epsilon), \quad \mbox{for $\phi_\epsilon \in H^1(\Omega_\epsilon)$.}
$$
With some abuse of notation we will identify all different realizations of this operator and we will write them all as $A_\epsilon$.

Under these assumptions, we write (\ref{nbc}) in an abstract form as $A_{\epsilon} u_\epsilon = h_\epsilon(u_\epsilon)$, where $h_\epsilon : H^1(\Omega_\epsilon) \to H^{-\alpha}(\Omega_\epsilon)$ with $\half <\alpha < 1$ is defined by
\begin{equation} 
\label{hepsilon}
\langle h_\epsilon(u_\epsilon), \phi_\epsilon\rangle = \frac{1}{\epsilon}\int_{\omega_\epsilon} f(x,u_\epsilon)\phi_\epsilon - \int_{\partial\Omega_\epsilon}g(x,u_\epsilon) \phi_\epsilon,  \quad \mbox{for $\phi_\epsilon \in H^\alpha(\Omega_\epsilon)$}.
\end{equation}

In particular, $u_\epsilon$ is a solution of (\ref{nbc}) if and only if $u_\epsilon$ satisfies $u_\epsilon=A_\epsilon^{-1}h_\epsilon(u_\epsilon)$, that is, $u_\epsilon$ is a fixed point of the nonlinear map $A_\epsilon^{-1}\circ h_\epsilon:H^1(\Omega_\epsilon)\to H^1(\Omega_\epsilon)$.

Similarly, the solutions of the limiting problem (\ref{nbc_limite_gamma_F}) can be written as fixed points of the map $A^{-1}_{0}\circ h_{0}: H^1(\Omega)\to H^1(\Omega)$, where $A_{0}:D(A_{0}) \subset L^2(\Omega) \to L^2(\Omega)$ is the linear operator defined by $A_{0} u = - \Delta u + u$ with domain $D(A_{0}) = \{u \in H^2(\Omega):  \frac{\partial u}{\partial n} = 0 \mbox{\ in\ } \partial\Omega\}$ and $h_{0}: H^1(\Omega) \to H^{-\alpha}(\Omega)$ with $\half <\alpha < 1$ is defined by 
\begin{equation} 
\label{hzero}
\langle h_{0}(u), \phi\rangle = \int_{\partial \Omega} \beta f(x,u)\phi -
\int_{\partial\Omega}\gamma g(x,u) \phi,  \quad \mbox{for $\phi \in H^\alpha(\Omega)$}.
\end{equation}


\subsection{Bessel Potentials spaces and trace operators}
\label{sobolev-trace}

The results obtained here are valid for the  more general case  of domain perturbation, which was considered in \cite{AB0}. For this we will assume that the hypothesis {\bf (P)} holds. Observe that $\Omega_\epsilon\cap U_i$ is a Lipschitz deformation of the fixed domain $Q_{N-1}\times (-1,0)$ and the Lipschitz norm of the transformation is uniformly controlled in $\epsilon$ as $\epsilon\to 0$. This fact  allow to obtain uniform estimates in Bessel Pontentials spaces and of constant of trace operators. 

We denote by $H^{s,p}(V)$, where $s \leq 1$, $1<p<\infty$ and $V$ is an open set or $\R^N$, the Bessel Potentials spaces. We note that $H^{0,p}(V)=L^{p}(V)$ and $H^{s,2}(V)=H^{s}(V)$. Let us start with the following result

\begin{lemma}
\label{mudancauniformeBessel}
Let $S$ and $S_\epsilon$, $0<\epsilon \leq \epsilon_0$, be families of bounded domains in $\R^N$. Assume there exists a family of Lipschitz one-to-one mappings $\Psi_\epsilon$ from  $S$  onto $S_\epsilon$  such that its inverse $\Psi_\epsilon^{-1}$ is Lipschitz, $\| D\Psi_\epsilon\|_{\mathcal{L}(L^{\infty}(S))} \leq K$ and $\|D\Psi_\epsilon^{-1}\|_{\mathcal{L}(L^{\infty}(S_\epsilon))} \leq K$, for some $K>0$ independent of $\epsilon$. Then $u \in H^{s,p}(S_\epsilon)$ iff $u\circ \Psi_\epsilon \in  H^{s,p}(S)$ and there exist constants $C$, $D>0$ independent of $\epsilon$ such that
\begin{equation}
\label{besselmaps}
C\|u\circ \Psi_\epsilon\|_{H^{s,p}(S)} \leq \|u\|_{H^{s,p}(S_\epsilon)} \leq D \|u\circ \Psi_\epsilon\|_{H^{s,p}(S)}.
\end{equation}
\end{lemma}
\noindent {\bf Proof. } Let $X^0= L^p(S)$, $X^0_\epsilon= L^p(S_\epsilon)$, $X^1= W^{1,p}(S)$ and $X^1_\epsilon= W^{1,p}(S_\epsilon)$. 
The spaces  $H^{s,p}(S)$ and $H^{s,p}(S_\epsilon)$ are the $X^\theta$ and $X^\theta_{\epsilon}$ complex interpolation spaces between $X^0$ and $X^1$, $X^0_\epsilon$ and $X^1_\epsilon$, respectively, which are defined below. 

Denote by $H(X^0, X^1)$ all continuous bounded functions $F:S\longrightarrow X^0$ which are holomorphic in the interior of $S=\{z \in \mathbb{C} \; : \; 0\leq \mbox{Re}(z)\leq 1\}$ such that the following norm is finite
$$
{\ds\|F\| = \max \left\{ \sup_{y}\|F( iy)\|_{X^0}, \sup_{y}\|F(1+ iy)\|_{X^1} \right\}}.
$$
For $0\leq \theta\leq 1$, the interpolation space $X^\theta$ is defined by $X^\theta=\{ F(\theta): F \in H(X^0, X^1)\}$ with norm 
$
\| f\|_\theta =\inf  \{ \| F\|: F(\theta)=f\}.
$
Similarly we define $X^\theta_{\epsilon}$. For more details see \cite{Amann1}.

Using the mappings $\Psi_\epsilon$ and $\Psi_\epsilon^{-1}$ we get, for $\epsilon$ fixed, that $u \in H^{s,p}(S_\epsilon)$ iff $u \circ \Psi_\epsilon \in  H^{s,p}(S)$, with $s\leq 1$. Now, using the norms of the interpolation spaces, $H^{s,p}(S)$ and $H^{s,p}(S_\epsilon)$, and applying Lemma 4.1 in \cite{AB0}, we  get  (\ref{besselmaps}). \quad $\blacksquare$ 

\vspace{0.3cm}

\begin{remark}
\label{boundedness-for-Phi_i_eps}
As an important example of mappings satisfying the hypotheses of Lemma \ref{mudancauniformeBessel} we mention the family of maps 
$$
\Phi_{i,\epsilon}:Q_N^-\to U_i\cap \Omega_\epsilon,
$$
where we denote by $Q_N^-=Q_{N-1}\times (-1,0)\subset Q_N$ and we assume hypothesis {\bf (P)}  is satisfied.
\end{remark}

With this lemma we can prove a result about the continuity of the trace operator for Bessel Potentials spaces. 

\begin{proposition} 
\label{remarkimersaotraco} 
Let $\Omega_\epsilon$, $0\leq\epsilon\leq\epsilon_0$, be a family of domains satisfying {\bf (P)}. Then the constant of the trace operator $H^{s,p}(\Omega_\epsilon)\to L^q(\partial\Omega_\epsilon)$, for $1<p<\infty$, $\frac{1}{p}<s\leq 1$ and $1\leq q\leq \frac{p(N-1)}{N-sp}$, is uniformly bounded in $\epsilon$, that is, there exists $C>0$ independent of $\epsilon$ such that  
$$
\|u_\epsilon\|_{L^q(\partial\Omega_\epsilon)}\leq C\|u_\epsilon\|_{H^{s,p}(\Omega_\epsilon)}, \quad \mbox{for all $u_\epsilon\in H^{s,p}(\Omega_\epsilon)$}. 
$$
\end{proposition}
\noindent {\bf Proof. } Observe that $\|u_\epsilon\|_{L^q(\partial\Omega_\epsilon)}\leq \sum_{i=1}^m\|u_\epsilon\|_{L^q(\partial\Omega_\epsilon\cap U_i)}$. Using Lemma \ref{mudancauniformeBessel} and the continuity of the trace operator for a fixed domain (see \cite{necas}), we have 
\begin{eqnarray*}
\|u_\epsilon\|_{L^q(\partial\Omega_\epsilon \cap U_i)} & \leq  & C\|u_\epsilon \circ \Phi_{i,\epsilon}\|_{L^q(Q_{N-1})} \leq \displaystyle \tilde{C}\|u_\epsilon \circ \Phi_{i,\epsilon}\|_{H^{s,p}(Q_{N}^-)} \\ & \leq&  \displaystyle \hat{C}\|u_\epsilon\|_{H^{s,p}(\Omega_\epsilon \cap U_i)}\leq \displaystyle \hat{C}\|u_\epsilon\|_{H^{s,p}(\Omega_\epsilon)}.  \quad  \blacksquare
\end{eqnarray*}

Moreover, we have

\begin{lemma} 
\label{uniforminterpolation}
Let $\Omega_\epsilon$, $0\leq\epsilon\leq\epsilon_0$, be a family of domains satisfying {\bf (P)}. Then the constant of the interpolation between $H^{s}(\Omega_\epsilon)$, $H^{1}(\Omega_\epsilon)$ and $L^2(\Omega_\epsilon)$, for $0<s<1$,  is uniformly bounded in  $\epsilon$. 
\end{lemma}
\noindent {\bf Proof. } In fact, using Lemma \ref{mudancauniformeBessel} and interpolation for a fixed domain,  we obtain
\begin{eqnarray*}
\|v\|_{H^{s}(\Omega_\epsilon \cap U_i)}  \leq  D\|v\circ \Phi_{i,\epsilon}\|_{H^{s}(Q_N^-)}  &\leq & \bar{C} \|v\circ \Phi_{i,\epsilon}\|^{1-s}_{L^{2}(Q_{N}^-)} \|v\circ \Phi_{i,\epsilon}\|^{s}_{H^{1}(Q_{N}^-)} \\ &\leq&  \tilde{C}\|v\|^{1-s}_{L^{2}(\Omega_\epsilon \cap U_i)} \ \|v\|^{s}_{H^{1}(\Omega_\epsilon \cap U_i)}.
\end{eqnarray*}
Hence, the result follows. \quad $\blacksquare$

\vspace{0.3cm}


\subsection{Concentrated integrals}
\label{results}

Now, since the result on trace operators has been established, we analyze  how the concentrated integrals converge to boundary integrals, as $\epsilon \to 0$. These convergence results will be needed to study the behavior of the nonlinear parts given by (\ref{hepsilon}) and (\ref{hzero}), and consequently to analyze the limit of the solutions of (\ref{nbc}). Initially, we have 

\begin{lemma} 
\label{lemconcentrat1} 
Assume hypotheses {\bf (P)} and {\bf(B)} are satisfied. Suppose that $v_{\epsilon}\in H^{s}(\Omega_{\epsilon})$ with $\frac{1}{2}<s\leq 1$ and $1\leq q\leq \frac{2(N-1)}{N-2s}$. Then, for small $\epsilon_{0}$, there exists a constant $C>0$ independent of $\epsilon$ and $v_{\epsilon}$ such that for any $0<\epsilon\leq\epsilon_{0}$, we have
\begin{equation}
\label{ineqconcentrat1}
\frac{1}{\epsilon}\int_{\omega_{\epsilon}}\left|v_{\epsilon}\right|^{q}  \leq C \left\|v_{\epsilon}\right\|^{q}_{H^{s}(\Omega_{\epsilon})}.
\end{equation}
\end{lemma}
\noindent {\bf Proof. } Consider the finite  cover $\{U_i\}_{i=0}^m$ such that $\overline{\Omega}_{\epsilon} \subset  \cup_{i=0}^m U_i\equiv U$. Note that  there exists $M>0$ independent of $\epsilon$ such that
\begin{eqnarray*}
\frac{1}{\epsilon}\!\int_{\omega_{\epsilon}\cap U_i}\!\!\!\left|v_{\epsilon}\right|^q \!\!\!\!
& = &\!\!\!\!    \frac{1}{\epsilon}\int_0^\epsilon \int_{Q_{N-1}}\!\! \left|v_\epsilon \left(\Phi_{i} \left( x',\sigma {\ds\frac{\rho_{i,\epsilon}(x')}{\epsilon}}\right) \right) \right|^q (J_N \Phi_i) \left(x^\prime, \sigma {\ds\frac{\rho_{i,\epsilon}(x')}{\epsilon}}\right) \frac{\rho_{i,\epsilon}(x')}{\epsilon} dx^\prime d\sigma \\
\!\!\!&  \leq &\!\!\!  \frac{M}{\epsilon}\int_0^\epsilon \int_{Q_{N-1}} \left|v_\epsilon \left(\Phi_{i} \left(x',\sigma {\ds\frac{\rho_{i,\epsilon}(x')}{\epsilon}}\right) \right) \right|^q J_{N-1}\Phi_i \left(x^\prime, \sigma {\ds\frac{\rho_{i,\epsilon}(x')}{\epsilon}}\right) dx^\prime d\sigma\\
\!\!\!& = &\!\!\! \frac{M}{\epsilon}\int^{\epsilon}_{0} \int_{\Gamma_{\sigma,\epsilon}}\left|v_{\epsilon}\right|^q d\sigma=\frac{M}{\epsilon}\int^{\epsilon}_{0} \left\|v_{\epsilon}\right\|^q_{L^{q}(\Gamma_{\sigma,\epsilon})} d\sigma,
\end{eqnarray*}
where $\Gamma_{\sigma,\epsilon}= \Phi_{i} \left( \left\{\left(x',\sigma {\ds\frac{\rho_{i,\epsilon}(x')}{\epsilon}}\right) \;:\; x' \in Q_{N-1}\right\} \right)$, $\sigma \in [0,\epsilon)$. 

Let $\Omega_{\sigma,\epsilon}= \Phi_{i} \left( \left\{\left(x',s {\ds\frac{\rho_{i,\epsilon}(x')}{\epsilon}}\right) \;:\; x' \in Q_{N-1}, \ s \in [0, \sigma)\right\} \right)$. Since $\Omega_{\sigma,\epsilon}$ satisfy {\bf(P)} then, using Proposition \ref{remarkimersaotraco} with $p=2$, we have that there exists $C>0$ independent of $\sigma$ and $\epsilon$ such that $\left\|v_{\epsilon}\right\|_{L^{q}(\Gamma_{\sigma,\epsilon})}\leq C \left\|v_{\epsilon}\right\|_{H^{s}(\Omega_{\sigma,\epsilon})}.$ Using the Sobolev embedding $H^{s}(\Omega_{\epsilon}) \hookrightarrow H^{s}(\Omega_{\sigma,\epsilon})$, we get (\ref{ineqconcentrat1}). \quad $\blacksquare$

\vspace{0.2cm}

With respect to the function $\beta$, given by Definition \ref{definition-of-beta}, we have the following

\begin{lemma}
\label{convconcent}
Assume hypotheses {\bf (P)} and {\bf(B)} are satisfied. Then, for any functions $h, \varphi \in H^{s}(U)$ with $\frac{1}{2} < s \leq 1$, we get 
\begin{eqnarray} 
\label{limiteexplicito}
\lim_{\epsilon \to 0}\frac{1}{\epsilon}\int_{\omega_{\epsilon}} h \, \varphi \, d\xi = \int_{\partial \Omega} \beta \, h \, \varphi \, dS.
\end{eqnarray}
\end{lemma}
\noindent {\bf Proof. } Consider the finite  cover $\{U_i\}_{i=0}^m$ such that $\overline{\Omega}_{\epsilon} \subset  \cup_{i=0}^m U_i\equiv U$. 
Initially, let $h$ and $\varphi$ be smooth functions defined in $\overline{U}$. Note that, using Definition \ref{definition-of-beta}, we get 
\begin{equation} 
\label{integralfronteira}
\begin{array}{lll}
\displaystyle \int_{\partial \Omega \cap U_i} \beta \, h \, \varphi \, dS & = &  \displaystyle \int_{Q_{N-1}} \beta(\Phi_{i}(x',0)) \, h(\Phi_{i}(x',0)) \, \varphi(\Phi_{i}(x',0)) J_{N-1}\Phi_{i}(x',0) \, dx'\\
& = & \displaystyle \int_{Q_{N-1}} \, h(\Phi_{i}(x',0)) \, \varphi(\Phi_{i}(x',0)) (J_{N}\Phi_{i})(x',0) \, \tilde{\beta}_i(x') \, dx',
\end{array}
\end{equation}

\begin{equation}
\label{integralconcentraweps}
\begin{array}{lcl}
{\ds{ \frac{1}{\epsilon}\int_{\omega_{\epsilon}\cap U_i}}}\!\!\! h \, \varphi \, d\xi  
\!\!\!\!& = & \!\!{\ds{\frac{1}{\epsilon} \int_{Q_{N\!-\!1}}\!\! \int^{\rho_{i,\epsilon}(x')}_{0}}}\!\!\! h(\Phi_{i}(x',x_{N})) \, \varphi(\Phi_{i}(x',x_{N})) J_{N}\Phi_{i}(x',x_{N}) \, dx_N dx' \\
& = & \!\!\!\!\displaystyle \int_{Q_{N\!-\!1}}\!\! \int^{1}_{0} h(\Phi_{i}(x',\rho_{i,\epsilon}(x')z)) \, \varphi(\Phi_{i}(x',\rho_{i,\epsilon}(x')z)) J_{N}\Phi_{i}(x', \rho_{i,\epsilon}(x')z) \frac{\rho_{i,\epsilon}(x')}{\epsilon} \, dz dx',
\end{array}
\end{equation}
where we take $x_N=\rho_{i,\epsilon}(x')z$. Since $\rho_{i,\epsilon}(x')z \to 0$ as $\epsilon\to 0$, uniformly for $(x',z) \in  Q_{N-1}\times[0,1]$, we have 
\begin{equation}
\label{c1}
\begin{array}{lll}
\left| \displaystyle \int_{Q_{N-1}} \int^{1}_{0}  \frac{\rho_{i,\epsilon}(x')}{\epsilon} [h(\Phi_{i}(x', \rho_{i,\epsilon}(x')z)) \, \varphi(\Phi_{i}(x', \rho_{i,\epsilon}(x')z)) J_{N}\Phi_{i}(x', \rho_{i,\epsilon}(x')z) \right. \\
\displaystyle \left. - h(\Phi_{i}(x',0)) \, \varphi(\Phi_{i}(x',0)) (J_{N}\Phi_{i})(x',0)]  \, dz dx'\right| \to 0, \quad \mbox{as $\epsilon \to 0$}.
\end{array}
\end{equation}

Now, from hypothesis {\bf(B)} we have
\begin{equation}
\label{c2}
\lim_{\epsilon \to 0} \int_{Q_{N-1}}\!\!\!\!\!h(\Phi_{i}(x',0)) \, \varphi(\Phi_{i}(x',0)) (J_{N}\Phi_{i})(x',0)  \,  {\ds\frac{\rho_{i,\epsilon}(x^\prime)}{\epsilon}} \, dx^\prime =\displaystyle \int_{Q_{N-1}}\!\!\!\!\! h(\Phi_{i}(x',0)) \, \varphi(\Phi_{i}(x',0)) (J_{N}\Phi_{i})(x',0)\, \tilde{\beta}_i(x') \, \, dx'.
\end{equation}

From \eqref{integralfronteira}, \eqref{integralconcentraweps}, (\ref{c1}) and (\ref{c2}), we obtain $\left|\frac{1}{\epsilon}\int_{\omega_{\epsilon}} h \, \varphi \, d\xi - \int_{\partial \Omega} \beta \, h \, \varphi \, dS\right| \to 0$ as $\epsilon \to 0$. 
Consequently, the proof of equality (\ref{limiteexplicito}) follows from density arguments and the continuity of the trace operator.  \quad $\blacksquare$

\vspace{0.2cm}

As a consequence of this result, we get

\begin{corollary}
\label{inde2}
Assume hypotheses {\bf (P)} and {\bf(B)} are satisfied. Then the function $\beta$ is independent of the parametrization chosen and therefore it is unique. 
\end{corollary}
\noindent {\bf Proof. } Suppose that $\beta$ depends on the parametrization. Then there will exist $\beta$ and  $\bar{\beta}$, both satisfying Lemma \ref{convconcent} with $\varphi \equiv 1$. Hence, by uniqueness of the limit, we have  
$$
\int_{\partial\Omega} \beta \; h \; dS = \int_{\partial\Omega} \bar{\beta}\; h \; dS ,  \quad \mbox{for all } h \in C_0^\infty(\R^N).
$$
This implies that $\beta =\bar{\beta}$ almost everywhere in $\partial\Omega$. \quad $\blacksquare$


\section{Upper semicontinuity of solutions}
\label{upper}

In this section we will provide a proof of Theorem \ref{th-main}, that is, we will show the upper semicontinuity of the family of solutions of (\ref{nbc}) and (\ref{nbc_limite_gamma_F}) in $H^{1}(\Omega_\epsilon)$. In particular, we will obtain that the limit problem of (\ref{nbc}) is given by (\ref{nbc_limite_gamma_F}). For this, we will keep the notation of the previous sections and, in terms of the nonlinearities, taking into account Remark \ref{rm-cut-off}, we will assume that $f$ and $g$ satisfy (\ref{hip-f}). Since $\Omega \subset \Omega_\epsilon$, $\Omega_\epsilon$ is an exterior perturbation of $\Omega$, we consider the restriction operator $R_\Omega: H^1(\Omega_\epsilon) \to H^1(\Omega)$ given by $R_\Omega(u)=u_{|\Omega}$. 

Let us start with the following result about $E$-convergence of a sequence of functions in $H^s(\Omega_\epsilon)$.

\begin{lemma} 
\label{convdeEepsilontilde}
Assume {\bf (P)} holds. Let $u_\epsilon \in H^1(\Omega_\epsilon)$ such that $\|u_\epsilon\|_{H^1(\Omega_\epsilon)} \leq M $, for some constant $M>0$ independent of $\epsilon$, and $ R_\Omega(u_{\epsilon}) \stackrel{\epsilon \to 0}{-{\hspace{-2mm}}\rightharpoonup}  u$ in $H^1(\Omega)$, then $\|u_{\epsilon} -  E_{\epsilon} u\|_{H^s(\Omega_{\epsilon})}\to 0$  as $\epsilon \to 0$, for $0<s<1$.
\end{lemma}
\noindent {\bf Proof. } Using Lemma \ref{uniforminterpolation} we get
$$
\|u_{\epsilon} -E_{\epsilon} u\|_{H^s(\Omega_{\epsilon})} \leq   M \|u_{\epsilon} -E_{\epsilon} u\|^{1-s}_{L^2(\Omega_{\epsilon})} \|u_{\epsilon} -E_{\epsilon} u\|^{s}_{H^1(\Omega_{\epsilon})}  \leq  M (\|R_\Omega u_{\epsilon} - u\|_{L^2(\Omega)} + \|u_{\epsilon} -E_{\epsilon} u\|_{L^2(\omega_{\epsilon})} )^{1-s}\|u_{\epsilon} -E_{\epsilon} u\|^{s}_{H^1(\Omega_{\epsilon})}.
$$
Hence, using that $\|u_\epsilon\|_{H^1(\Omega_\epsilon)}$ and $\left\|E_\epsilon u\right\|_{H^{1}(\Omega_{\epsilon})}$ are bounded uniformly in $\epsilon$, $R_\Omega u_{\epsilon} \to u$ in $L^2(\Omega)$ and $|\omega_\epsilon| \to 0$ when $\epsilon \to 0$, H\"older Inequality and uniform embedding of $H^1(\Omega_\epsilon)$ in $L^q(\Omega_\epsilon)$, for $1\leq q < \frac{2N}{N-2} $ (see Proposition 4.2 in \cite{AB0}), the result follows. \quad $\blacksquare$

\vspace{0.2cm}

Now, we prove a result that will be used to obtain the boundedness of solutions of (\ref{nbc}).

\begin{lemma}
\label{limitacaodanorma}
Assume {\bf (P)} and {\bf (B)} are  satisfied. Let $\{u_\epsilon\}$ be a sequence in  $H^{1}(\Omega_\epsilon)$ and let $z_\epsilon$ be given by $z_{\epsilon} = A_\epsilon^{-1} h_\epsilon(u_\epsilon)$. Then $\{z_{\epsilon}\}$ is a bounded sequence in $H^1(\Omega_\epsilon)$.
\end{lemma}
\noindent {\bf Proof. } Recall that saying that  $z_{\epsilon} = A_\epsilon^{-1} h_\epsilon(u_\epsilon)$ is equivalent to saying that $z_\epsilon$ is the weak solution of 
\begin{equation}
\label{parte1}
\left\{
\begin{array}{ll}
-\Delta z_\epsilon + z_\epsilon=\displaystyle \frac{1}{\epsilon} \mathcal{X}_{\omega_{\epsilon}} f(x,u_\epsilon), & \hbox{in} \ \Omega_\epsilon \\
\frac{\ts\partial z_\epsilon}{\ts\partial n}+ g(x,u_\epsilon) =0, & \hbox{on} \ \partial \Omega_\epsilon.
\end{array} \right.
\end{equation}
Multiplying the equation (\ref{parte1}) by $z_\epsilon$, integrating by parts and using (\ref{hip-f}), Lemma \ref{lemconcentrat1} and Proposition \ref{remarkimersaotraco}, we get 
$$
\|z_\epsilon\|^2_{H^1(\Omega_\epsilon)}  = \displaystyle \frac{1}{\epsilon} \int_{\omega_\epsilon} f(x,u_\epsilon)z_\epsilon -
\int_{\partial\Omega_\epsilon} g(x,u_\epsilon)z_\epsilon \leq \displaystyle C_1\left\|z_{\epsilon}\right\|_{H^{1}(\Omega_{\epsilon})}+C\left\|z_{\epsilon}\right\|_{L^{1}(\partial \Omega_{\epsilon})} \leq \displaystyle C_1\left\|z_{\epsilon}\right\|_{H^{1}(\Omega_{\epsilon})}+C_2\left\|z_{\epsilon}\right\|_{H^{1}(\Omega_{\epsilon})}.
$$
Hence, there exists $K>0$ independent of $\epsilon$ such that $\|z_\epsilon\|_{H^1(\Omega_\epsilon)}\leq K$. \quad $\blacksquare$

\vspace{0.2cm}

In order to obtain the upper semicontinuity of the family of solutions of (\ref{nbc}) and (\ref{nbc_limite_gamma_F}), we study the behavior of the nonlinearities $h_\epsilon$, $0\leq \epsilon \leq \epsilon_0$, defined by (\ref{hepsilon}) and (\ref{hzero}). Here, we also proved a more general result,  which the behavior of  potential $V_\epsilon$,  $0\leq \epsilon \leq \epsilon_0$, defined in $\Omega_\epsilon$  is analyzed when $\epsilon \to 0$. This result will be used in Section \ref{conclusion}.  

\begin{proposition} 
\label{convpartenaolinearna}
Assume {\bf (P)} and {\bf(B)} are satisfied. Let $\{u_\epsilon\}$  be a bounded sequence in $H^1(\Omega_\epsilon)$ such that  $R_{\Omega}(u_\epsilon)\stackrel{\epsilon \to 0}{-{\hspace{-2mm}}\rightharpoonup} u$ in $H^1(\Omega)$.\\
{\bf (i) } If $\{z_\epsilon\}$ is a bounded sequence in $H^1(\Omega_\epsilon)$ such that  $R_{\Omega}(z_\epsilon)\stackrel{\epsilon \to 0}{-{\hspace{-2mm}}\rightharpoonup}  z$  in $H^1(\Omega)$,  then 
$$
\begin{array}{l}
\displaystyle \frac{1}{\epsilon}\int_{\omega_\epsilon} f(x,u_\epsilon)z_\epsilon\to
\int_{\partial \Omega}\beta f(x,u)z, \quad \mbox{as $\epsilon \to 0$}. 
\end{array}
$$
{\bf (ii)} If $V_\epsilon$  is a potential defined in $\Omega_\epsilon$  such that $\|V_\epsilon\|_{L^\infty(\omega_{\epsilon})} \leq K$, for some $K>0$ independent of $\epsilon$, and $V_\epsilon \circ \psi_{i} \to V_0\circ \psi_{i} $ in $L^2(Q_{N-1})$, for all $i=1,2,\ldots, m$,   then
$$
\frac{1}{\epsilon} \int_{\omega_{\epsilon}}V_{\epsilon}u_{\epsilon}v  \to  \int_{\partial \Omega}\beta V_0 uv, \quad \mbox{as $\epsilon \to 0$, \quad $\forall v\in H^1(U)$}.
$$
\end{proposition}

\noindent {\bf Proof. } {\bf (i)} Using (\ref{hip-f}), Cauchy-Schwarz and Lemma \ref{lemconcentrat1}, we have
$$
\begin{array}{lll}
\displaystyle \left|\frac{1}{\epsilon}\int_{\omega_\epsilon} f(x,u_\epsilon)z_\epsilon-
\int_{\partial\Omega} \beta f(x,u)z \right|\\
\\
\displaystyle \leq \frac{1}{\epsilon}\int_{\omega_\epsilon}|f(x,u_\epsilon)-f(x,E_\epsilon u)||z_\epsilon| +\frac{1}{\epsilon}\int_{\omega_\epsilon} |f(x,E_\epsilon u)||z_\epsilon-E_\epsilon z|  \displaystyle +\left|\frac{1}{\epsilon}\int_{\omega_\epsilon} f(x,E_\epsilon u)E_\epsilon z-\int_{\partial\Omega}\beta f(x,u)z\right|  \\
\\
\displaystyle \leq \frac{1}{\epsilon}\int_{\omega_\epsilon}|\partial_{u}f(x,\theta_{\epsilon}u_\epsilon+(1-\theta_{\epsilon})E_\epsilon u)|\left|u_\epsilon-E_\epsilon u\right||z_\epsilon| +\frac{C}{\epsilon}\int_{\omega_\epsilon} |z_\epsilon-E_\epsilon z| +\left|\frac{1}{\epsilon}\int_{\omega_\epsilon} f(x,E_\epsilon u)E_\epsilon z-\int_{\partial\Omega}\beta f(x,u)z\right|\\
\\
\displaystyle \leq C \left( \frac{1}{\epsilon}\int_{\omega_\epsilon}\left|u_\epsilon- E_\epsilon u\right|^2 \right)^{\frac{1}{2}} \left(\frac{1}{\epsilon}\int_{\omega_\epsilon}|z_\epsilon|^2\right)^{\frac{1}{2}}  +\frac{C}{\epsilon}\int_{\omega_\epsilon} |z_\epsilon- E_\epsilon z| +\left|\frac{1}{\epsilon}\int_{\omega_\epsilon} f(x,E_\epsilon u)E_\epsilon z-\int_{\partial\Omega}\beta f(x,u)z\right|  \\
\\
\displaystyle \leq C\left\|u_{\epsilon}-E_\epsilon u\right\|_{H^{s}(\Omega_{\epsilon})} \left\|z_{\epsilon}\right\|_{H^{1}(\Omega_{\epsilon})} + C\left\|z_{\epsilon}-E_\epsilon z\right\|_{H^{s}(\Omega_{\epsilon})}+ \left|\frac{1}{\epsilon}\int_{\omega_\epsilon} f(x,E_\epsilon u)E_\epsilon z-\int_{\partial\Omega}\beta f(x,u)z\right| \to 0, \quad \mbox{as $\epsilon \to 0$,}
\end{array}
$$
for some $0\leq \theta_{\epsilon}(x) \leq 1$, $x\in \Omega_\epsilon$, where we use  Lemma \ref{convconcent} to prove that the last term goes to $0$. Now, choosing $s\in \R$ such that $\frac{1}{2}<s<1$ and using Lemma \ref{convdeEepsilontilde}, we get $\left\|u_{\epsilon}-E_\epsilon u\right\|_{H^{s}(\Omega_{\epsilon})} \to 0$ and $\left\|z_{\epsilon}-E_\epsilon z\right\|_{H^{s}(\Omega_{\epsilon})} \to 0$, as $\epsilon \to 0$. \quad 

\vspace{0.2cm}

\noindent {\bf (ii)}  Initially, using that $V_\epsilon$ is uniformly bounded in $L^\infty(\omega_\epsilon)$, Cauchy-Schwarz and Lemma \ref{lemconcentrat1}, we get
$$
\begin{array}{lll}
\displaystyle \left|\frac{1}{\epsilon}\int_{\omega_{\epsilon}}V_{\epsilon}u_{\epsilon}v -\int_{\partial\Omega}\beta V_0 uv\right| & \leq &  \displaystyle \frac{1}{\epsilon}\int_{\omega_{\epsilon}}\left|V_{\epsilon}\right|\left|u_{\epsilon}-E_{\epsilon} u\right|\left|v\right| +\left|\frac{1}{\epsilon}\int_{\omega_{\epsilon}}V_{\epsilon}E_{\epsilon} uv- \int_{\partial\Omega}\beta V_0 uv\right| \\
& \leq & \displaystyle K \left\|u_{\epsilon}-E_{\epsilon} u\right\|_{H^{s}(\Omega_{\epsilon})}\left\|v\right\|_{H^{1}(\Omega_{\epsilon})} +\left|\frac{1}{\epsilon}\int_{\omega_{\epsilon}}V_{\epsilon}E_{\epsilon} uv- \int_{\partial\Omega}\beta V_0 uv\right|.
\end{array}
$$
Choosing $s\in \R$ such that $\frac{1}{2}<s<1$ and using Lemma \ref{convdeEepsilontilde}, we get $\left\|u_{\epsilon}-E_{\epsilon} u\right\|_{H^{s}(\Omega_{\epsilon})} \to 0$ as $\epsilon \to 0$. 

Consider the finite  cover $\{U_i\}_{i=0}^m$ such that $\overline{\Omega}_{\epsilon} \subset  \cup_{i=0}^m U_i\equiv U$. For each $i$, we denote by  $\Phi_{i\rho_\epsilon}(x^\prime,z)=\Phi_i(x^\prime,\rho_{i,\epsilon}(x^\prime) z)$. We have, 
\begin{eqnarray*}
& & \left|\frac{1}{\epsilon}\int_{\omega_{\epsilon}\cap U_i}V_{\epsilon}E_{\epsilon}uv- \int_{\partial\Omega\cap U_i}\beta V_0 uv\right| \\
& = &\displaystyle \left|  \int_{\!\!Q_{N-1}}\!\int_0^{1} \!V_\epsilon(\Phi_{i\rho_\epsilon}(x^\prime,z))E_\epsilon u (\Phi_{i\rho_\epsilon}(x^\prime,z)) v(\Phi_{i\rho_\epsilon}(x^\prime,z))J_N\Phi_i(x^\prime,\rho_{i,\epsilon}(x^\prime) z)\frac{\rho_{i,\epsilon} (x^\prime)}{\epsilon} dz dx^\prime \right. \\
& - & \left. \int_{Q_{N-1}} \ V_0(\psi_{i}(x^\prime))\ u(\psi_{i}(x^\prime))\ v(\psi_{i}(x^\prime)) \tilde{\beta}_i(x^\prime)(J_{N}\Phi_{i})(x^\prime,0)dx^\prime \right|\\
& \leq &  \left|\int_{Q_{N-1} }\!\!\! \frac{\rho_{i,\epsilon}(x')}{\epsilon}\! \int_0^{1}\!\!V_\epsilon(\Phi_{i\rho_\epsilon}(x^\prime, z))E_\epsilon u (\Phi_{i\rho_\epsilon}(x^\prime,z)) v(\Phi_{i\rho_\epsilon}(x^\prime,z)) [J_N\Phi_i(x^\prime,\rho_{i,\epsilon}(x^\prime) z)\! -\! (J_N\Phi_i)(x^\prime,0) ]dz dx^\prime \right| \\
&+ & \left|  \int_{Q_{N-1}} \!\!\frac{\rho_{i,\epsilon}(x^\prime)}{\epsilon}\! \int_0^{1}\!V_\epsilon(\Phi_{i\rho_\epsilon}(x^\prime, z))E_\epsilon u(\Phi_{i\rho_\epsilon}(x^\prime, z))  [v(\Phi_{i\rho_\epsilon}(x^\prime, z))-v(\Phi_i(x^\prime,0)) ] (J_N\Phi_i)(x^\prime,0) dz dx^\prime\right| \\
&+&
\left| \int_{Q_{N-1}}\!\! \frac{\rho_{i,\epsilon}(x^\prime)}{\epsilon}\! \int_0^{1}\!V_\epsilon(\Phi_{i\rho_\epsilon}(x^\prime, z)) [E_\epsilon u( \Phi_{i\rho_\epsilon}(x^\prime, z))-  E_\epsilon u(\Phi_i(x^\prime,0)) ] v(\Phi_i(x^\prime,0))  (J_N\Phi_i)(x^\prime,0)  dz dx^\prime \right| \\
&+&
\left| \int_{Q_{N-1}}\!\! \frac{\rho_{i,\epsilon}(x^\prime)}{\epsilon}\! \int_0^{1}\!V_\epsilon(\Phi_{i\rho_\epsilon}(x^\prime, z)) [E_\epsilon u( \Phi_{i\rho_\epsilon}(x^\prime, z))-  E_\epsilon u(\Phi_i(x^\prime,0)) ] v(\Phi_i(x^\prime,0))  (J_N\Phi_i)(x^\prime,0)  dz dx^\prime \right| 
\end{eqnarray*}
\begin{eqnarray*}
&+& 
\left| \int_{Q_{N-1}}\!\! \frac{\rho_{i,\epsilon}(x^\prime)}{\epsilon}\! \int_0^{1}\![V_\epsilon(\Phi_{i\rho_\epsilon}(x^\prime, z))- V_\epsilon(\Phi_i(x^\prime,0))] E_\epsilon u (\Phi_i(x^\prime,0)) v(\Phi_i(x^\prime,0)) (J_N\Phi_i)(x^\prime,0) dz dx^\prime \right|\\
&+&
 \left|\int_{Q_{N-1}}\! \frac{\rho_{i,\epsilon}(x^\prime)}{\epsilon} \int_0^{1}\![V_\epsilon(\Phi_i(x^\prime,0))-V_0(\Phi_i(x^\prime,0))] E_\epsilon u(\Phi_i(x^\prime,0)) v(\Phi_i(x^\prime,0)) (J_N\Phi_i)(x^\prime,0) dz dx^\prime \right|\\
&+& 
\left|\int_{Q_{N-1}} \left(\frac{\rho_{i,\epsilon}(x^\prime)}{\epsilon}  -\tilde{\beta}_i(x^\prime) \right) V_0(\psi_i(x^\prime)) u(\psi_i(x^\prime)) v(\psi_i(x^\prime))  (J_{N}\Phi_i)(x^\prime,0) dx^\prime \right| = I_1+I_2+I_3+I_4+I_5+I_6.
\end{eqnarray*}

Using  Lemma 4.2 in \cite{AB0} we get that  $ v(\Phi_i(x^\prime,\rho_{i,\epsilon}(x^\prime) z)) \to v(\psi_i(x^\prime))$, $E_\epsilon u(\Phi_i(x^\prime,\rho_{i,\epsilon}(x^\prime) z)) \to E_\epsilon u(\psi_i(x^\prime))$ and   $ |V_\epsilon (\Phi_i(x^\prime,\rho_{i,\epsilon}(x^\prime) z)) - V_\epsilon(\psi_i(x^\prime))| \to 0$ in $L^q(Q_{N-1})$, $1\leq q< \frac{2N-2}{N-2}$, uniformly in $z  \in [0,1]$.  Since $V_\epsilon$ is uniformly bounded in $L^\infty(\omega_\epsilon)$,  then $ |V_\epsilon (\Phi_i(x^\prime,\rho_{i,\epsilon}(x^\prime) z)) - V_\epsilon(\psi_i(x^\prime))| \to 0$ in $L^r(Q_{N-1})$, $1\leq r<\infty$. Hence, $I_2$, $I_3$ and $I_4$ converge to zero when $\epsilon \to 0$. 
Moreover,  $J_N \Phi_i$ is continuous and $\rho_{i,\epsilon}(x^\prime)
 $ goes to zero, then  $I_1$ converges to zero when $\epsilon \to 0$. 
Since $V_\epsilon \circ \psi_{i} \to V_0\circ \psi_{i} $ in $L^2(Q_{N-1})$, then $I_5$ goes to zero. 
Finally, using hypothesis {\bf(B)}, we obtain that $I_6$ goes to zero. \quad $\blacksquare$

\begin{proposition}
\label{convpartenaolinearnaoexterior}
Assume {\bf (P)}, {\bf (B)} and {\bf (G)} are satisfied. Let $\{u_\epsilon\}$ and $\{z_\epsilon\}$ be bounded sequences in $H^1(\Omega_\epsilon)$ such that $R_{\Omega}(u_\epsilon) \stackrel{\epsilon \to 0}{-{\hspace{-2mm}}\rightharpoonup} u$ and $R_{\Omega}(z_\epsilon) \stackrel{\epsilon \to 0}{-{\hspace{-2mm}}\rightharpoonup}  z$ both in $H^1(\Omega)$. Then 
$\langle h_\epsilon (u_\epsilon), z_\epsilon\rangle \to \langle h_0(u), z\rangle$, as $\epsilon \to 0$.
\end{proposition}
\noindent {\bf Proof. } By Proposition 5.1 in \cite{AB0}, 
$$
\displaystyle \int_{\partial\Omega_\epsilon} g(x,u_\epsilon)z_\epsilon \to
\int_{\partial\Omega} \gamma g(x,u)z , \quad \mbox{as $\epsilon \to 0$}.
$$
Hence, using  Proposition \ref{convpartenaolinearna}, we get $\langle h_\epsilon (u_\epsilon), z_\epsilon\rangle \to \langle h_0(u), z\rangle$, as $\epsilon \to 0$. \quad $\blacksquare$

\vspace{0.2cm}

So, we are in conditions to prove the upper semicontinuity of the family of solutions of (\ref{nbc}) and (\ref{nbc_limite_gamma_F}).

\vspace{0.2cm}

\par\medskip\noindent {\bf Proof of Theorem \ref{th-main}:} If $\{u_\epsilon^*\}$, $0<\epsilon \leq \epsilon_0$, is a family of solution of (\ref{nbc}) satisfying $\| u_\epsilon^*\|_{L^\infty(\Omega_\epsilon)}\leq R$, for some constant $R>0$ independent of $\epsilon$, then $u_\epsilon^*=A^{-1}_{\epsilon}h_{\epsilon}(u_\epsilon^*)$ and, by Lemma \ref{limitacaodanorma}, we have that $\{u_\epsilon^*\}$ is a bounded sequence in $H^{1}(\Omega_{\epsilon})$. Therefore,  by Lemma 4.3  in \cite{AB0},  there exist a subsequence $\{u_{\epsilon_k}^{*}\}$ and $u_0^* \in H^{1}(\Omega)$ such that 
$$
R_{\Omega}(u^{*}_{\epsilon_k})  \stackrel{\epsilon \to 0}{-{\hspace{-2mm}}\rightharpoonup} u^{*}_0 \quad \mbox{in $H^1(\Omega)$} \qquad \mbox{and} \qquad u^{*}_{\epsilon_k}\dwto u^{*}_0.
$$

Let us show that $u^{*}_0$ is a solution of (\ref{nbc_limite_gamma_F}), that is, $u^{*}_0=A^{-1}_{0}h_{0}(u^{*}_0)$. For this, notice that since $u^{*}_{\epsilon_k}\dwto u^{*}_0$, we have 
$$
(u^{*}_{\epsilon_k},E_{\epsilon_k}v)_{H^1(\Omega_{\epsilon_k})}\to (u^{*}_0,v)_{H^1(\Omega)}, \quad \forall v\in H^1(\Omega).
$$ 
But, using Proposition \ref{convpartenaolinearnaoexterior} we get
$$
(u^{*}_{\epsilon_k}, E_{\epsilon_k}v)_{H^1(\Omega_{\epsilon_k})} = \langle h_{\epsilon_k}(u^{*}_{\epsilon_k}), E_{\epsilon_k} v\rangle \to \langle h_{0}(u^{*}_0), v \rangle = (A^{-1}_{0}h_{0}(u^{*}_0),v)_{H^1(\Omega)}, \quad \forall v\in H^1(\Omega).
$$
Hence, $u^{*}_0= A^{-1}_{0}h_{0}(u^{*}_0)$. 

Now, we prove that $u^*_{\epsilon_k} {\stackrel{E}{\longrightarrow}} u^*_0$. In order to do this, we prove the convergence of the norms  $\|u^{*}_{\epsilon_k}\|_{H^1(\Omega_{\epsilon_k})} \to \|u^{*}_0\|_{H^1(\Omega)}$. Using again  Proposition \ref{convpartenaolinearnaoexterior}, we have
$$
\|u^{*}_{\epsilon_k}\|^2_{H^1(\Omega_{\epsilon_k})} = \langle h_{\epsilon_k}(u^{*}_{\epsilon_k}), u^{*}_{\epsilon_k}\rangle  \to \langle h_{0}(u^{*}_0), u^{*}_0\rangle=\|u^{*}_0\|_{H^1(\Omega)}^2.
$$
 The convergence of the norms and the weak $E$-convergence of the sequence imply, by Proposition 3.2 in \cite{AB0},  that $u^{*}_{\epsilon_k} \dto u^{*}_0$. \quad $\blacksquare$


\section{Final conclusion}
\label{conclusion}

With the results obtained in this work and proceeding analogously to Corollary 5.3 in \cite{AB0}, we can prove the lower semicontinuity of the family of solutions of (\ref{nbc}) and (\ref{nbc_limite_gamma_F}) in $H^1(\Omega_\epsilon)$, in the case where the solution of the limit problem (\ref{nbc_limite_gamma_F}) is hyperbolic. Also, we  leave out  the proof that for any hyperbolic solution of the limit problem (\ref{nbc_limite_gamma_F}), there exists one and only one solution of (\ref{nbc}) in its neighborhood, since it follows from Proposition \ref{convpartenaolinearna} (ii) and Proposition 5.4 in \cite{AB0}  and by similar arguments to Proposition 5.5 in \cite{AB0}. Moreover, using Proposition \ref{convpartenaolinearna} (ii), we can also prove that if $\{u_\epsilon^*\}$ is a sequence of solutions of (\ref{nbc}) which converge to $u_0^*$, a solution of (\ref{nbc_limite_gamma_F}), then the eigenvalues and eigenfunctions of the linearization of (\ref{nbc}) around $u_\epsilon^*$ converge to the eigenvalues and eigenfunctions of the linearization of (\ref{nbc_limite_gamma_F}) around $u_0^*$.


\end{document}